\documentclass[12pt,fleqn,leqno]{article}
\usepackage{graphicx,color,enumerate,amssymb}
\usepackage{amsmath}
\usepackage{amssymb}
\usepackage{epsfig}
\usepackage{array}
\usepackage{float}

\makeatletter
\@addtoreset{equation}{section}
\renewcommand\thetable{\@arabic\c@table}
\renewcommand\thefigure{\@arabic\c@figure}
\long\def\@makecaption#1#2{%
  \vskip\abovecaptionskip
  \begin{center}%
  \sbox\@tempboxa{#1: #2}%
  \ifdim \wd\@tempboxa >\hsize
    #1: #2\par
  \else
    \global \@minipagefalse
    \hb@xt@\hsize{\hfil\box\@tempboxa\hfil}%
  \fi
  \end{center}%
  \vskip\belowcaptionskip}
\def\N{{\rm I\kern-.15em N}}
\def\R{{\rm I\kern-.2em R}}
\def\Z{{\rm Z\kern-.26em Z}}
\makeatother
\newtheorem{thm}{Theorem}[section]

\newtheorem{prop}[thm]{Proposition}

\newcommand{\be}{\begin{eqnarray}}
\newcommand{\ee}{\end{eqnarray}}
\newcommand{\bq}{\begin{eqnarray*}}
\newcommand{\eq}{\end{eqnarray*}}

\newcommand{\verk}{\stackrel{{\cal D}}{\longrightarrow}}

\newcommand{\bewend}{\hspace*{2mm}\rule{3mm}{3mm}}

\newcommand{\RR}{\mathbb{R}}
\newcommand{\PP}{\mathbb{P}}

\newcommand{\BE}{\mathbb{E}}
\newcommand{\BV}{\mathbb{V}}

\newcommand{\vertk}{\stackrel{{\cal D}}{\longrightarrow}}

\begin{document}
\begin{center}
{ \LARGE\sc On the consistency of the  spacings test for multivariate uniformity}\\
 \vspace*{0.5cm} {\large\sc Norbert Henze}\\ \vspace*{0.5cm}
{\it Institute of Stochastics, Karlsruhe Institute of Technology, Karlsruhe, Germany} \\
\end{center}

{\small {\bf Abstract.}{ We give a simple conceptual proof of the consistency of a test for multivariate uniformity in a bounded set $K \subset \RR^d$
that is based on the maximal spacing generated by i.i.d. points $X_1,\ldots,X_n$ in $K$, i.e., the volume of the largest convex set of a given shape that is contained in $K$ and avoids
each of these points. Since asymptotic results for the case $d >1$ are only availabe under uniformity, a key element of the proof is a suitable coupling.}\\

\vspace*{0.3cm}

{\small {\it Keywords.}
Multivariate spacings; test for multivariate uniformity; consistency; extreme-value distribution} \\
{2010 AMS classification:} 62H15, 62G20.

\section{Introduction and main result}\label{sec_1}
Whereas there is a plethora of procedures for testing the hypothesis that a random sample comes from the uniform distribution on the unit interval,
the problem of testing for uniformity of multivariate observations has hitherto been studied far less fully. This paper is not concerned
with giving an overview over the available literature (see e.g. \cite{BCV2006}, \cite{EHJ2016}, \cite{Jupp2008}, \cite{PEWI2013}), but to turn attention to the
maximum spacings test studied in \cite{BCPL2012}. To be specific, let $K$ be a bounded set in $\RR^d$, $d \ge 1$, where $|K| =1$ and $|\partial K| =0$ and
$|\cdot|$ denotes  Lebesgue measure. Moreover, let $A\subset \RR^d$ be a fixed bounded convex set with nonempty interior and $|A|=1$. If $X_1,X_2, \ldots $ are
independent and identically distributed (i.i.d.) random vectors on some common probability space $(\Omega,{\cal A},\PP)$
taking values in $K$, the maximum spacing
 (with respect to the reference set $A$) generated by $X_1,\ldots,X_n$
is defined
as
\[
\Delta_n = \sup \{r : \exists\, x \textrm{ with } x+rA \subset K \setminus \{X_1,\ldots,X_n\}\}.
\]
Letting $V_n = \Delta_n^d$ denote the $d$-dimensional volume (Lebesgue measure) of this 'maximum gap' defined by $X_1,\ldots,X_n$, \cite{Jan1987} proved the following:
If the distribution of $X_1$ is uniform over $K$, then
\begin{equation}\label{janson}
\lim_{n\to \infty} \PP\left( nV_n - \log n - (d-1)\log \log n - \beta \le t \right) = G(t), \quad t \in \RR,
\end{equation}
where $G(t) = \exp(-\exp(-t))$ is the  distribution function of the extreme value distribution of Gumbel, and $\beta$ is some constant that depends only on the boundary of $A$ (and is zero if $A$ is a cube).
\cite{Jan1987} also proved
\begin{eqnarray*}
\liminf_{n \to \infty} \frac{nV_n- \log n}{\log \log n } = d-1 \quad \PP\textrm{-a.s.},\\
\limsup_{n \to \infty} \frac{nV_n- \log n}{\log \log n } = d+1 \quad \PP\textrm{-a.s.},
\end{eqnarray*}
thus giving values conjectured in \cite{Deh1983} and generalizing earlier results to non-cubical gaps. Notice that the above limits imply
\begin{equation}\label{aslimit}
\lim_{n\to \infty} \frac{n V_n}{\log n} = 1 \quad \PP\textrm{-a.s.}
\end{equation}

In the case $d=1$ and $K=[0,1]$, result (\ref{janson}) is due to L. Weiss, see \cite{Weiss1960}. This paper has been largely forgotten, since
it is referenced neither in \cite{Jan1987} nor in \cite{BCPL2012}. The latter paper suggests to use $V_n$ as a statistic for testing the hypothesis $H_0$ that
$X_1$ has a uniform distribution over $K$. Using (\ref{janson}) and denoting by $g_{1-\alpha}$ the ($1-\alpha$)-quantile of $G$, an asymptotic level-$\alpha$-test
rejects $H_0$ if $V_n > c_{n,\alpha}$,
where
\[
c_{n,\alpha} = \frac{g_{1-\alpha} +  \log n + (d-1) \log \log n + \beta}{n}.
\]
In the univariate case, this test has been proposed in \cite{Weiss1960}, and \cite{Weiss1960} also proves its consistency against general alternatives.
The authors of \cite{BCPL2012} prove the consistency of the test based on $V_n$ (see Theorem 1 of \cite{BCPL2012}). This proof, however, hinges on a heuristic argument (see line 4
 of the left-hand column of page 270 of [1]), and the method of proof does not cover the case of testing for uniformity on the surface of a sphere or, more generally, on
 a lower-dimensional differentiable manifold of $\RR^d$.

 It is the purpose of this paper to give a conceptual, simple proof of the  consistency of the maximal spacings test against general alternatives.
Our main result is as follows.

\vspace*{3mm}

\begin{thm} \label{mainth}
Suppose that $X_1$ has a Lebesgue density $f$. If there is some $\varepsilon >0$ and some sphere $S \subset K$ so that $f(x) \le 1- \varepsilon$ for each $x \in S$, then
\begin{equation}\label{convone}
\lim_{n \to \infty} \PP(V_n > c_{n,\alpha}) = 1.
\end{equation}
\end{thm}
Thus, the test based on $V_n$ is consistent against each such alternative.

\section{Proof of Theorem \ref{mainth}}
To prove Theorem \ref{mainth}, a key observation is that -- under a suitable condition -- a random sample size does not change the limit  behaviour (\ref{janson}).
To this end, let $O_\PP(a_n)$ denote a random variable that is bounded in probability after division by $a_n$, where $(a_n)$ is a sequence of positive real numbers.
Furthermore, write $\verk$ for convergence in distribution.

\begin{prop}\label{proprs} Let $(k_n)$ be a sequence of integers satisfing
\begin{equation}\label{constab}
0< a \le \frac{k_n}{n} \le b < \infty, \quad n \ge 1,
\end{equation}
for some constants $a,b$. Furthermore, let $(L_n)$ be a sequence of integer-valued
random variables defined on $(\Omega,{\cal A},\PP)$ so that $L_n = k_n + O_\PP(\sqrt{n})$. Then
\[
k_n V_{L_n} - \log k_n - (d-1) \log \log k_n - \beta \vertk G \quad \textrm{as } n \to \infty,
\]
where {\rm{(}}by an abuse of notation{\rm{)}} the random variable $G$ has a Gumbel distribution.
\end{prop}

In other words, (\ref{janson}) continues to hold if the fixed sample size is replaced by a random one, provided that both sample sizes do not differ too much.
\vspace*{3mm}

  \noindent {\sc Proof} of Proposition \ref{proprs}. From (\ref{janson}) we have
\[
L_n V_{L_n} - \log L_n - \log \log L_n - \beta \vertk G \quad \textrm{as } n \to \infty.
\]
Now, (\ref{constab}) and $L_n = k_n + O_\PP(\sqrt{n})$ yield
\[
\log L_n = \log k_n + o_\PP(1), \quad \log \log L_n = \log \log k_n + o_\PP(1).
\]
Furthermore, we have
\[
(L_n -k_n)V_{L_n} = \frac{L_n-k_n}{\sqrt{n}} \cdot \frac{L_n V_{L_n}}{\log L_n} \cdot \frac{k_n}{L_n} \cdot \frac{\sqrt{n}\log L_n}{k_n}.
\]
The first factor on the right-hand side is $O_\PP(1)$, and the second converges to 1 almost surely in view of (\ref{aslimit}). Since the third and the
last factor are $1+ o_\PP(1)$ and $o_\PP(1)$, repectively, the assertion follows from Sluzki's lemma.  \bewend

\vspace*{2mm}

\noindent {\sc Proof} of Theorem \ref{mainth}: We will prove this theorem by means of a suitable coupling
that 'mediates' between the uniform and the alternative distribution. To be specific, let $S$ and $\varepsilon$ be as in the statement of the theorem. Writing
${\bf 1}\{A\}$ for the indicator function of an event $A$, let
\[
N_n = \sum_{j=1}^n {\bf 1}\{X_j \in S\}
\]
be the number of points that fall into $S$. The distribution of $N_n$ is binomial with parameters $n$ and $p$, where
\[
p = \PP(X_1 \in S) = \int_S f(x)\, \textrm{d}x.
\]
We will assume $p>0$ since otherwise $\liminf_{n\to \infty} V_n >0$ $\PP$-almost surely, from which (\ref{convone}) follows. Let
$Z_1,\ldots,Z_{N_n}$ denote those $X_j$ that fall into $S$, and write $W_{N_n}$ for the volume of the largest spacing generated by $Z_1,\ldots,Z_{N_n}$ within $S$.
Then
\begin{equation}\label{ungl}
W_{N_n} \le V_n.
\end{equation}
Without loss of generality let the underlying probability space $(\Omega,{\cal A},\PP)$ be rich enough to carry a sequence
$(Y_j,U_j)_{j \ge 1}$ of i.i.d. random vectors, which are independent of $X_1,X_2, \ldots$, each having a uniform distribution on $S \times [0,1-\varepsilon]$.
Notice that $Y_j$ and $U_j$ are independent, and that the distributions of $Y_j$ and $U_j$ are uniform over $S$ and $[0,1-\varepsilon]$, respectively.
The crucial point is that the conditional distribution of $Y_1$ given that $U_1 \le f(Y_1)$ equals the distribution of $Z_1$ (this fact is also known
as the acceptance-rejection method in connection with the generation of random numbers from a distribution with a bounded density).
We now observe $(Y_1,U_1), (Y_2,U_2), \ldots $ sequentially. If $U_i \le f(Y_i)$ for the $j$th time, we denote the corresponding $Y_i$ by $\widetilde{Z}_j$.
Let
\[
L_n = \min \bigg{\{} k \ge 1: \sum_{i=1}^k {\bf 1}\{U_i \le f(Y_i)\} = N_n \bigg{\}}
\]
be number of samples from $(Y_1,U_1), (Y_2,U_2), \ldots $ needed to obtain $N_n$ points
$\widetilde{Z}_1, \ldots, \widetilde{Z}_{N_n}$. We then have
\begin{equation}\label{edist}
(\widetilde{Z}_1, \ldots, \widetilde{Z}_{N_n}) \ \stackrel{{\cal D}}{=} \ (Z_1,\ldots,Z_{N_n}),
\end{equation}
where $\stackrel{{\cal D}}{=}$ stands for equality in distribution.
Letting $\widetilde{W}_{N_n}$ denote the volume of the maximal spacing generated by $\widetilde{Z}_1, \ldots, \widetilde{Z}_{N_n}$ within $S$, and writing
$\widetilde{V}_{L_n}$ for the volume of the largest spacing generated by $Y_1,\ldots,Y_{L_n}$ within $S$, the fact that $L_n \ge N_n$ and (\ref{edist}), (\ref{ungl}) imply
$
\widetilde{V}_{L_n} \le \widetilde{W}_{N_n} \stackrel{{\cal D}}{=} W_{N_n} \le V_n.
$
Thus, (\ref{convone}) would follow if we can show
\begin{equation}\label{convone1}
\lim_{n \to \infty} \PP(\widetilde{V}_{L_n} > c_{n,\alpha}) = 1.
\end{equation}
To this end, we need some information on the random variable $L_n$ in order to be able to apply Proposition \ref{proprs}. Since $L_n$ models the waiting time until
$N_n$ 'successes', i.e., 'cases $U_i \le f(Y_i)$',  have been obtained, we have
\[
L_n \ \stackrel{{\cal D}}{=} \ \sum_{j=1}^{N_n} (1+\xi_j),
\]
where $\xi_1,\xi_2, \ldots $ is a sequence of i.i.d. random variables, which is also independent of $N_n$, and
$\xi_1$ has a geometric distribution with parameter $\kappa$ (say), where
\[
\kappa  =   \PP(U_1 \le f(Y_1)) = \BE\big{[} \PP(U_1 \le f(Y_1)|Y_1)\big{]}  = \BE \bigg{[} \frac{f(Y_1)}{1-\varepsilon} \bigg{]}
 =  \frac{p}{(1-\varepsilon)|S|}.
\]
It follows that
\[
\BE(L_n) = \BE(N_n) \cdot \BE(1+\xi_1) = np \cdot \frac{1}{\kappa} = n(1-\varepsilon)|S|
\]
and (using $\BV(\xi_1) = (1-\kappa)/\kappa^2$)
\begin{eqnarray*}
\BV(L_n) & = & \BE[\BV(L_n|N_n)] + \BV(\BE[L_n|N_n]) \\
& = & np \cdot \frac{1-\kappa}{\kappa^2} + \frac{1}{\kappa^2} \cdot np(1-p)\\
& = & n \cdot \frac{(1-\varepsilon)^2|S|^2}{p} \left(2-p-\frac{p}{(1-\varepsilon)|S|}\right).
\end{eqnarray*}
From Chebyshev's inequality, we thus have $L_n = k_n + O_\PP(\sqrt{n})$ as $n \to \infty$,
where $k_n = \lfloor n (1-\varepsilon)|S| \rfloor $ and $\lfloor \cdot \rfloor$ denotes the floor function. Now, Proposition \ref{proprs} yields (notice that we have to multiply the volume of the maximum spacing
by $1/|S|$)
\[
k_n \frac{\widetilde{V}_{L_n}}{|S|} - \log k_n - (d-1) \log \log k_n - \beta \vertk G \quad \textrm{as } n \to \infty.
\]
The definition of $k_n$ and Sluzki's lemma then entail
\[
n(1-\varepsilon) \widetilde{V}_{L_n} - \log n - (d-1) \log \log n - \log((1-\varepsilon)|S|) - \beta \vertk G \quad \textrm{as } n \to \infty.
\]
Abbreviating the expression to the left of '$\vertk$' by $V^*_{L_n}$, the definition of $c_{n,\alpha}$ gives
\[
\PP(\widetilde{V}_{L_n} > c_{n,\alpha}) = \PP\left(V^*_{L_n} > - \varepsilon \, \log n \cdot (1+ a_n)\right)
\]
for some sequence $(a_n)$ converging to $0$, and (\ref{convone1}) follows. \bewend

\section{Concluding remarks}
\begin{enumerate}
\item [1)] The proof of consistency reveals that the test based on $V_n$ is consistent against any alternative distribution that gives probability zero
to some nonempty open subset $O$ of $K$, since then $\liminf_{n\to \infty} V_n >0$ $\PP$-almost surely.
\item[2)]  Result (\ref{janson}) remains true for spherical spacings on a sphere, and more generally for geodesic balls on any compact $C^2$-Riemannian
manifold, see the second remark on p. 276 of \cite{Jan1987}. Our proof of consistency is general enough to carry over almost literally to cover also
these cases.
\item[3)] There is an analogue to the largest multivariate spacing, which is the largest nearest neighbour distance. Letting $\|\cdot \|$   denote the Euclidean norm in
$\RR^d$, let $D_{n,i} = \min_{j\neq i}\|X_i-X_j\|$ denote the nearest neighbour distance of $X_i$ to the remaining points. To avoid boundary effects which  may dominate
in higher dimensions (see \cite{Sttj1986}, \cite{DeHe1989}, \cite{DeHe1990}), let
\[
M_n = \max_{i=1,\ldots,n} \min\left(D_{n,i}, \textrm{dist}(X_i, \partial K)\right),
\]
where dist$(X_i,\partial K)$ is the distance of $X_i$ to the boundary of $K$.
Thus, $M_n$ is the radius of the largest sphere contained in $K$ that has one of the points as center and avoids all other points. Letting $\overline{V}_n$ denote the
volume of the sphere with radius $M_n$, \cite{Hen1983} proved
\[
n\overline{V}_n  - \log n \vertk G \quad \textrm{as } n \to \infty
\]
and showed consistency of a test for uniformity that rejects $H_0$ for large values of $\overline{V}_n$ against general alternatives (for a generalization to $r$th nearest neighbours, see \cite{Hen1982}).
 Moreover,
if $X_{n,1}, \ldots, X_{n,n}$, $n \ge 1$, is a triangular array of rowwise i.i.d. random vectors with density
\begin{equation}\label{contig}
f_n(x) = 1 + \frac{h(x)}{\log n}, \quad x \in K,
\end{equation}
where $h:K \rightarrow \RR$ is a nonnegative continuous function satisfying $\int_K h(x) \, \textrm{d} x =0$, then, under this
sequence of alternatives, we have
\[
\lim_{n\to \infty} \PP\left( n \overline{V}_n - \log n \le t\right) = G(t-C(h)), \quad t \in \RR,
\]
where
\[
C(h) = \log \int_K \exp(-h(x)) \, \textrm{d}x.
\]
Since $C(h) > 0$ if $h \not\equiv 0$, the test has positive asymptotic power against contiguous alternatives of the type (\ref{contig}).
We conjecture that the test based on $V_n$ shares this property.
\end{enumerate}




\vspace*{3mm}
\noindent Norbert Henze, Institute of Stochastics\\
Karlsruhe Institute of Technology (KIT)\\
 Englerstr. 2, D-76133 Karlsruhe\\
Norbert.Henze@kit.edu

\end{document}